\documentclass[11pt,reqno]{amsart}
\usepackage{amssymb,amscd,amsbsy,mathrsfs}
\renewcommand{\a}{\alpha}
\renewcommand{\b}{\beta}
\newcommand{\g}{\gamma}
\renewcommand{\d}{\delta}

\newcommand{\z}{\zeta}

\newcommand{\s}{\sigma}

\newcommand{\f}{\varphi}

\newcommand{\G}{\Gamma}
\newcommand{\D}{\Delta}

\newcommand{\h}{{\mathscr H}}

\newcommand{\X}{{\mathscr X}}

\newcommand{\T}{{\Bbb T}}

\newcommand{\Z}{{\Bbb Z}}

\newcommand{\0}{{\boldsymbol{0}}}

\newcommand{\bS}{{\boldsymbol S}}

\newcommand{\rf}[1]{(\ref{#1})}

\newcommand{\df}{\stackrel{\mathrm{def}}{=}}

\newcommand{\trace}{\operatorname{trace}}
\newcommand{\rank}{\operatorname{rank}}

\newcommand{\eeq}{\end{equation}}
\newcommand{\beq}{\begin{equation}}
\newcommand{\bay}{\begin{eqnarray}}
\newcommand{\ba}{\begin{align*}}
\newcommand{\ea}{\end{align*}}
\newcommand{\ey}{\end{eqnarray}}
\newcommand{\bey}{\begin{eqnarray*}}
\newcommand{\eey}{\end{eqnarray*}}

\newcommand{\be}{\infty}

\newcommand{\bl}{\blacksquare}
\newcommand{\ess}{\operatorname{ess}}

\newcommand{\Pf}{{\bf Proof. }}

\newcommand{\ov}{\overline}

\newtheorem{thm}{\hspace{\parindent}Theorem}[section]

\newtheorem{lem}[thm]{\hspace{\parindent}Lemma}

\usepackage{amsmath} 

\makeatletter
\def\upintkern@{\mkern-7mu\mathchoice{\mkern-3.5mu}{}{}{}}
\def\upintdots@{\mathchoice{\mkern-4mu\@cdots\mkern-4mu}%
 {{\cdotp}\mkern1.5mu{\cdotp}\mkern1.5mu{\cdotp}}%
 {{\cdotp}\mkern1mu{\cdotp}\mkern1mu{\cdotp}}%
 {{\cdotp}\mkern1mu{\cdotp}\mkern1mu{\cdotp}}}

\newcommand{\UpMultiIntegral}[1]{%
  \edef\ints@c{\noexpand\upintop
    \ifnum#1=\z@\noexpand\upintdots@\else\noexpand\upintkern@\fi
    \ifnum#1>\tw@\noexpand\upintop\noexpand\upintkern@\fi
    \ifnum#1>\thr@@\noexpand\upintop\noexpand\upintkern@\fi
    \noexpand\upintop
    \noexpand\ilimits@
  }%
  \futurelet\@let@token\ints@a
}
\makeatother

\DeclareFontFamily{OMX}{mdbch}{}
\DeclareFontShape{OMX}{mdbch}{m}{n}{ <->s * [0.8]  mdbchr7v }{}
\DeclareFontShape{OMX}{mdbch}{b}{n}{ <->s * [0.8]  mdbchb7v }{}
\DeclareFontShape{OMX}{mdbch}{bx}{n}{<->ssub * mdbch/b/n}{}

\DeclareSymbolFont{uplargesymbols}{OMX}{mdbch}{m}{n}
\SetSymbolFont{uplargesymbols}{bold}{OMX}{mdbch}{b}{n}
\DeclareMathSymbol{\upintop}{\mathop}{uplargesymbols}{82}
\DeclareMathSymbol{\upointop}{\mathop}{uplargesymbols}{"48}

\DeclareFontEncoding{MDB}{}{}
\DeclareFontFamily{MDB}{mdbch}{}
\DeclareFontShape{MDB}{mdbch}{m}{n}{ <->s * [0.8]  mdbchrmb }{}
\DeclareFontShape{MDB}{mdbch}{b}{n}{ <->s * [0.8]  mdbchbmb }{}
\DeclareFontShape{MDB}{mdbch}{bx}{n}{<->ssub * mdbch/b/n}{}
\DeclareFontSubstitution{MDB}{cmr}{m}{n}
\DeclareSymbolFont{mathdesignB}{MDB}{mdbch}{m}{n}%
\SetSymbolFont{mathdesignB}{bold}{MDB}{mdbch}{b}{n}%
\DeclareMathSymbol{\upintclockwise}{\mathop}{mathdesignB}{128}
\DeclareMathSymbol{\upointclockwise}{\mathop}{mathdesignB}{130}
\DeclareMathSymbol{\upointctrclockwise}{\mathop}{mathdesignB}{132}
\DeclareMathSymbol{\upoiint}{\mathop}{mathdesignB}{134}
\DeclareMathSymbol{\upoiiint}{\mathop}{mathdesignB}{136}

\makeatletter
\newcommand{\upint}{\DOTSI\upintop\ilimits@}
\newcommand{\upoint}{\DOTSI\upointop\ilimits@}
\makeatother

\theoremstyle{remark}

\newtheorem*{rem*}{Remark}

\newcommand\mB{\mathcal{B}}

\newcommand\mT{\mathcal{T}}

\newcommand{\diag}{\operatorname{diag}}

\begin{document}

\newcommand{\vse}{\vspace{.2in}}
\numberwithin{equation}{section}

\title{Multiple operator integrals, Haagerup and Haagerup-like tensor products, and operator ideals}
\author{A.B. Aleksandrov and V.V. Peller}
\thanks{The first author is partially supported by RFBR grant 14-01-00198; the second author is partially supported by NSF grant DMS 1300924}
\thanks{Corresponding author: V.V. Peller; email: peller@math.msu.edu}

\begin{abstract}
We study Schatten--von Neumann properties of multiple operator integrals with 
integrands in the Haagerup tensor product of $L^\be$ spaces. We obtain sharp, best possible estimates. This allowed us to obtain sharp Schatten--von Neumann estimates in the case of {\it Haagerup-like} tensor products.
\end{abstract}

\maketitle

{\bf
\footnotesize
\tableofcontents
\normalsize
}

\

%\begin{center}
%{\Large Contents}
%\end{center}
%
%\
%
%\begin{enumerate}
%\item[1.] Introduction \quad\dotfill \pageref{In}
%\item[2.] Preliminaries \quad\dotfill \pageref{prel}
%\item[oooooooooobnbgggggbbbbbb 	3.] Triple operator integrals
%\quad\dotfill \pageref{Stoi}
%\item[4.] Schatten--von Neumann properties of triple operator integrals  \quad\dotfill \pageref{SvNtoi}
%\item[5.] Haagerup like tensor products and triple operator integrals \quad\dotfill \pageref{Ttoi}
%\item[6.] When do the divided differences $\dg^{[1]}f$ and 
%$\dg^{[2]}f$ belong to Haagerup-like tensor products? \quad\dotfill \pageref{ddiff}
%\item[7.] Lipschitz type estimates in the case $1\le p\le2$ \quad\dotfill \pageref{ple2}
%\item[8.] No Lipschitz type estimates in the operator norm and in the 
%$\bS_p$ norm for \lb$p>2\,\,$! \quad\dotfill \pageref{Bp>2}
%\item[9.] Two counterexamples \quad\dotfill \pageref{2c}
%\item[10.] Points of local Lipschitzness \quad\dotfill \pageref{tlL}
%\item[11.] A sufficient condition for Lipschitz type estimates \quad\dotfill \pageref{dostu}
%\item[12.] Functions of noncommuting unitary operators \quad\dotfill \pageref{unit}
%\item[] References \quad\dotfill \pageref{bibl}
%\end{enumerate}
%
%\

\setcounter{section}{0}
\section{\bf Introduction}
\setcounter{equation}{0}
\label{In}

\medskip

Double operator integrals have been used in perturbation theory for 60 years. They were introduced by Yu.L. Daletskii and S.G. Krein in \cite{DK}, where the authors studied the problem of differentiability of operator functions and found a formula for operator derivatives in terms of double operator integrals. Later Birman and Solomyak developed in \cite{BS1}, \cite{BS2}, and \cite{BS3} the beautiful theory of double operator integrals and applied their theory to a variety of problems.
Since that time double operator integrals were successfully applied in perturbation theory on many occasions (see, e.g., \cite{Pe1}, \cite{Pe2}, \cite{AP2}, \cite{AP3}, \cite{APPS}, \cite{NP} and \cite{PoS}).

It turned out that various problems of perturbation theory lead to multiple operator integrals. Several mathematicians tried to give definitions of multiple operator integrals, see \cite{Pa}, \cite{St}. However, those definitions required very strong restrictions on the classes of functions that can be integrated. In \cite{Pe3} multiple operator integrals were defined for functions that belong to the (integral) projective tensor product of $L^\be$ spaces. 
Such multiple operator integrals were used in the problem of evaluating higher operator derivatives (see \cite{Pe2} and \cite{Pe3}), in estimating higher order operator differences (see \cite{AP1} and \cite{AP2}). They also appeared in a natural way when studying various trace formulae for functions of self-adjoint operators 
(see \cite{AP3}).

Later in \cite{JTT} multiple operator integrals were defined for Haagerup tensor products of $L^\be$ spaces, which is broader that the projective tensor product of $L^\be$ spaces. Nevertheless, it turns out that multiple operator integrals defined for the Haagerup tensor products of $L^\be$ spaces do not possess such natural Schatten--von Neumann classes as in the case of projective tensor products.

Triple operator integrals played an important role recently in problems of estimating functions of noncommuting self-adjoint and unitary operators, see 
\cite{ANP1}, \cite{ANP2}, \cite{ANP3}, \cite{AP4} and \cite{AP5}.
In those papers new tensor products of $L^\be$ spaces were introduced (so-called Haagerup-like tensor products of the first kind and and Haagerup-like tensor products of the second kind) and
triple operator integrals were defined for integrands that belong to such Haagerup-like tensor products. 

Note that in \cite{ANP1}, \cite{ANP2}, \cite{ANP3}, \cite{AP4} and \cite{AP5} it was important to obtain Schatten--von Neumann estimates of triple operator integrals with integrands belonging to Haagerup and Haagerup-like tensor products of $L^\be$ spaces.

We also mention here the recent survey article \cite{Pe4} on multiple operator integrals in perturbation theory.

Note that the anonymous referee of \cite{ANP3} observed that the proof of the theorem on Schatten--von Neumann properties of triple operator integrals with integrands in the Haagerup tensor product of $L^\be$ spaces can be simplified and the simplification also allows one to extend the result of \cite{ANP3} to a broader range of Schatten--von Neumann ideals.  

In this paper we obtain most general results on Schatten--von Neuman properties of {\it arbitrary multiple operator integrals}
 with integrands belonging to Haagerup(-like) tensor products of $L^\be$ spaces and prove that these results are best possible. 

\

\section{\bf An introduction to multiple operator integrals}
\setcounter{equation}{0}
\label{prel}

\

In this section we recall the definition of triple operator integrals in the case when the integrand belongs to the projective (or integral projective) tensor product of $L^\be$ spaces. Then we proceed to multiple operator integrals with integrands in the Haagerup tensor product of $L^\be$.

Let us first consider, for simplicity, the case of triple operator integrals.

Let $E_1$, $E_2$, and $E_3$ be spectral measures on Hilbert space and let $T$ and $R$ be bounded linear operators on Hilbert space. Triple operator integrals are expressions of the following form:
\bay
\label{troi}
\int\limits_{\X_1}\int\limits_{\X_2}\int\limits_{\X_3} 
\Psi(x_1,x_2,x_3)\,dE_1(x_1)T\,dE_2(x_2)R\,dE_3(x_3).
\ey
Such integrals make sense under certain assumptions on $\Psi$, $T$, and $R$. The function $\Psi$ is called the {\it integrand} of the triple operator integral.

Recall that the {\it projective tensor product} 
$L^\be(E_1)\hat\otimes L^\be(E_2)\hat\otimes L^\be(E_3)$ can be defined as the class of function $\Psi$ of the form
\bay
\label{pred}
\Psi(x_1,x_2,x_3)=\sum_n\f_n(x_1)\psi_n(x_2)\chi_n(x_3)
\ey
such that
\bay
\label{norma}
\sum_n\|\f_n\|_{L^\be(E_1)}\|\psi_n\|_{L^\be(E_2)}\|\chi_n\|_{L^\be(E_3)}<\be.
\ey
The norm $\|\Psi\|_{L^\be\hat\otimes L^\be\hat\otimes L^\be}$ of $\Psi$ is, by definition, the infimum of the left-hand side of \rf{norma} over all representations of the form \rf{pred}.

For $\Psi\in L^\be(E_1)\hat\otimes L^\be(E_2)\hat\otimes L^\be(E_3)$ of the form 
\rf{pred} the triple operator integral \rf{troi} was defined in \cite{Pe3} by
\begin{align}
\label{otoi}
\iiint \Psi(x_1,x_2,x_3)&\,dE_1(x_1)T\,dE_2(x_2)R\,dE_3(x_3)\nonumber\\
=&\sum_n\left(\int\f_n\,dE_1\right)T\left(\int\psi_n\,dE_2\right)R\left(\int\chi_n\,dE_3\right).
\end{align}
Clearly, \rf{norma} implies that the series on the right converges absolutely in the norm. The right-hand side of \rf{otoi} does not depend on the choice of a representation of the form \rf{pred}, see, e.g., \cite{ACDS} and  \cite{Pe4}. Clearly,
\bay
\label{ograntoi}
\left\|\iiint \Psi(x_1,x_2,x_3)\,dE_1(x_1)T\,dE_2(x_2)R\,dE_3(x_3)\right\|
\le\|\Psi\|_{L^\be\hat\otimes L^\be\hat\otimes L^\be}\|T\|\cdot\|R\|.
\ey
Note that for $\Psi\in L^\be(E_1)\hat\otimes L^\be(E_2)\hat\otimes L^\be(E_3)$, triple operator integrals have the following Schatten--von Neumann properties:
\bay
\label{bep}
T\in\mB(\h),\quad R\in\bS_p,\quad1\le p<\be,\quad\Longrightarrow\quad 
\iiint\Psi\,dE_1T\,dE_2R\,dE_3\in\bS_p
\ey
and
\bay
\label{pq}
T\in\bS_p,\!\quad\! R\in\bS_q,\quad \frac1p+\frac1q\le1,\!
\!\quad\Longrightarrow\!\quad 
\iiint\!\Psi dE_1TdE_2RdE_3\in\bS_r,\quad\!\frac1r=\frac1p+\frac1q.
\ey
Here $\bS_p$, $p>0$, stands for the Schatten--von Neumann ideal of operators on Hilbert space. We refer the reader to \cite{GK} for definitions and properties of such ideals.

Let us also mention that the definition of multiple operator integrals was extended in \cite{Pe3} for 
functions $\Psi$ that belong to the so-called {\it integral projective tensor product} of the corresponding $L^\be$ spaces. We refer the reader to \cite{Pe3} for more detail.

Practically the same definition can be given for multiple operator integrals of the form
$$
\underbrace{\int\cdots\int}_m\Psi(x_1,\cdots,x_m)
\,dE_1(x_1)T_1\,dE_2(x_2)T_2\cdots T_{m-1}\,dE_m(x_m).
$$
Such multiple operator integrals can be defined for functions $\Psi$ in the projective tensor product $L^\be(E_1)\hat\otimes L^\be(E_2)\hat\otimes\cdots\hat\otimes L^\be(E_m)$ as well as for $\Psi$ in the integral projective tensor
(see \cite{Pe3}). 

Moreover, with such multiple operators integrals analogs of \rf{ograntoi}, \rf{bep} and \rf{pq} hold. In particular, 
$$
\underbrace{\int\cdots\int}_m\Psi(x_1,\cdots,x_m)
\,dE_1(x_1)T_1\,dE_2(x_2)T_2\cdots T_{m-1}\,dE_m(x_m)\in\bS_r,
$$
whenever 
$$
\Psi\in L^\be(E_1)\hat\otimes L^\be(E_2)\hat\otimes\cdots\hat\otimes L^\be(E_m),
$$
$$
T_1\in\bS_{p_1},~T_2\in\bS_{p_2},~\cdots,~T_{m-1}\in\bS_{p_{m-1}}\quad\mbox{and}
\quad\frac1r\df\frac1{p_1}+\frac1{p_2}+\cdots+\frac1{p_{m-1}}\le1.
$$

We proceed now to the approach to multiple operator integrals based on the Haagerup tensor product of $L^\be$ spaces. We refer the reader to the book \cite{Pi} for detailed information about Haagerup tensor products.

Let us first define the {\it Haagerup tensor product} $L^\be(E_1)\!\otimes_{\rm h}\!L^\be(E_2)$
of two $L^\be$ spaces as the space of functions $\Psi$ of the form
\bay
\label{Htprolya2}
\Psi=\sum_j\a_j(x_1)\b_j(x_2),\quad\mbox{where}\quad
\{\a_j\}_{j\ge0}\in L_{E_1}^\be(\ell^2)\quad\mbox{and}\quad
\{\b_j\}_{j\ge0}\in L_{E_2}^\be(\ell^2),
\ey
i.e.,
$$
\|\{\a_j\}_{j\ge0}\|_{L^\be(\ell^2)}\df
E_1\mbox{-}\ess\sup\left(\sum_{j\ge0}|\a_j(x_1)|^2\right)^{1/2}<\be
$$
and
$$
\|\{\b_j\}_{j\ge0}\|_{L^\be(\ell^2)}\df
E_2\mbox{-}\ess\sup\left(\sum_{j\ge0}|\b_j(x_1)|^2\right)^{1/2}<\be.
$$
The norm $\|\Psi\|_{L^\be\!\otimes_{\rm h}\!L^\be}$ is, by definition the infimum of the products
$$
\|\{\a_j\}_{j\ge0}\|_{L^\be(\ell^2)}\|\{\b_j\}_{j\ge0}\|_{L^\be(\ell^2)}
$$
over all representations of $\Psi$ in the form of \rf{Htprolya2}.

Note that the space $L^\be(E_1)\!\otimes_{\rm h}\!L^\be(E_2)$ coincides with the space ${\frak M}(E_1,E_2)$ of Schur multipliers with respect to $E_1$ and $E_2$ (see \cite{Pe1} and \cite{AP6}). For $\Psi\in L^\be(E_1)\!\otimes_{\rm h}\!L^\be(E_2)$ the double operator integral $\iint\Psi\,dE_1T\,dE_2$ can be defined
by
$$
\iint\Psi(x_1,x_2)\,dE_1(x_1)T\,dE_2(x_2)=
\sum_j\left(\int\a_j\,dE_1\right)T\left(\int\b_j\,dE_2\right),
$$
the series on the right converges in the weak operator topology, its sum does not depend on the choice of a representation of  $\Psi$ the form \rf{Htprolya2} and determines a bounded linear transformer on the space $\mB(\h)$ of bounded linear operators on 
$\h$ (see \cite{BS1}, \cite{Pe1} and \cite{AP6}).

Let us proceed now to triple operator integrals with integrands in the Haagerup tensor product of $L^\be$ spaces.
We define the {\it Haagerup tensor product}
$L^\be(E_1)\!\otimes_{\rm h}\!L^\be(E_2)\!\otimes_{\rm h}\!L^\be(E_3)$ as the space of function $\Psi$ of the form
\bay
\label{htr}
\Psi(x_1,x_2,x_3)=\sum_{j,k\ge0}\a_j(x_1)\b_{jk}(x_2)\g_k(x_3),
\ey
where $\a_j$, $\b_{jk}$, and $\g_k$ are measurable functions such that
\bay
\label{ogr}
\{\a_j\}_{j\ge0}\in L_{E_1}^\be(\ell^2), \quad 
\{\b_{jk}\}_{j,k\ge0}\in L_{E_2}^\be({\mathcal B}),\quad\mbox{and}\quad
\{\g_k\}_{k\ge0}\in L_{E_3}^\be(\ell^2),
\ey
where ${\mathcal B}$ is the space of matrices that induce bounded linear operators on $\ell^2$ and this space is equipped with the operator norm. In other words,
$$
\|\{\b_{jk}\}_{j,k\ge0}\|_{L^\be({\mathcal B})}\df
E_2\mbox{-}\ess\sup\|\{\b_{jk}(x_2)\}_{j,k\ge0}\|_{{\mathcal B}}<\be.
$$
By the sum on the right-hand side of \rf{htr} we mean 
$$
\lim_{M,N\to\be}~\sum_{j=0}^N\sum_{k=0}^M\a_j(x_1)\b_{jk}(x_2)\g_k(x_3).
$$
Clearly, the limit exists almost everywhere.

The norm of $\Psi$ in 
$L^\be\!\otimes_{\rm h}\!L^\be\!\otimes_{\rm h}\!L^\be$ is, by definition, 
the infimum of
$$
\|\{\a_j\}_{j\ge0}\|_{L^\be(\ell^2)}\|\{\b_{jk}\}_{j,k\ge0}\|_{L^\be({\mB})}
\|\{\g_k\}_{k\ge0}\|_{L^\be(\ell^2)}
$$
over all representations of $\Psi$ of the form \rf{htr}.

It is well known that $L^\be\hat\otimes L^\be\hat\otimes L^\be\subset
L^\be\!\otimes_{\rm h}\!L^\be\!\otimes_{\rm h}\!L^\be$,
see e.g., \cite{ANP3} and \cite{Pe4}.

In \cite{JTT} multiple operator integrals were defined for functions in the Haagerup tensor product of $L^\be$ spaces. Let 
$\Psi\in L^\be\!\otimes_{\rm h}\!L^\be\!\otimes_{\rm h}\!L^\be$ and  
suppose that \rf{htr} and \rf{ogr} hold. The triple operator integral \rf{troi} is defined by
\begin{align}
\label{htraz}
\iiint\Psi(x_1,x_2,x_3)&\,dE_1(x_1)T\,dE_2(x_2)R\,dE_3(x_3)\nonumber\\[.2cm]
=&
\sum_{j,k\ge0}\left(\int\a_j\,dE_1\right)T\left(\int\b_{jk}\,dE_2\right)
R\left(\int\g_k\,dE_3\right)\nonumber\\[.2cm]
=&\lim_{M,N\to\be}~\sum_{j=0}^N\sum_{k=0}^M
\left(\int\a_j\,dE_1\right)T\left(\int\b_{jk}\,dE_2\right)
R\left(\int\g_k\,dE_3\right).
\end{align}

Then the series in \rf{htraz} converges in the weak operator topology, the sum of the series does not depend on the choice of a representation
\rf{htr}, it determines a bounded linear operator and
\bay
\label{ogra}
\left\|\iiint\Psi(x_1,x_2,x_3)\,dE_1(x_1)T\,dE_2(x_2)R\,dE_3(x_3)\right\|
\le\|\Psi\|_{L^\be\!\otimes_{\rm h}\!L^\be\!\otimes_{\rm h}\!L^\be}
\|T\|\cdot\|R\|
\ey
(see e.g., \cite{ANP3}).

For completeness, we prove \rf{ogra} in the next section.

To define Haagerup tensor products of $L^\be$ spaces and corresponding multiple operator integrals in the general case, to simplify the notation, we consider quadruple operator integrals. It will be clear how to generalize the construction to the general case of multiple operator integrals.

Let $E_1$, $E_2$, $E_3$ and $E_4$ be spectral measures on a Hilbert space $\h$.
The {\it Haagerup tensor product} $L^\be(E_1)\!\otimes_{\rm h}\!L^\be(E_2)\!\otimes_{\rm h}\!L^\be(E_3)\otimes_{\rm h}\!L^\be(E_4)$ is the class of functions $\Psi$ of the form
\bay
\label{Haag4}
\Psi(x_1,x_2,x_3,x_4)=\sum_{j,k,l}\a_j(x_1)\b_{jk}(x_2)\g_{kl}(x_3)\d_l(x_4),
\ey
where
$$
\{\a_j\}_{j\ge0}\in L_{E_1}^\be(\ell^2),~\; 
\{\b_{jk}\}_{j,k\ge0}\in L_{E_2}^\be({\mathcal B}),~\;
\{\g_{kl}\}_{k,l\ge0}\in L_{E_3}^\be({\mathcal B}),
~\;\mbox{and}~\;
\{\d_l\}_{l\ge0}\in L_{E_4}^\be(\ell^2).
$$
The norm of $\Psi$ in $L^\be(E_1)\!\otimes_{\rm h}\!L^\be(E_2)\!\otimes_{\rm h}\!L^\be(E_3)\otimes_{\rm h}\!L^\be(E_4)$ is defined in the same way as in the case of three variables.

Quadruple operator integrals with integrands in the Haagerup tensor products of $L^\be$ spaces are defined by
\begin{align*}
\iiiint\!\Psi(x_1,x_2,x_3,x_4)dE_1(x_1)T_1dE_2(x_2)T_2dE_3(x_3)T_3dE_4(x_4)
\!=\!\sum_{j,k,l}A_jT_1B_{jk}T_2\G_{kl}T_3\D_l,
\end{align*}
where
$$
A_j\df\int\a_j\,dE_1,\quad
B_{jk}\df\int\b_{jk}\,dE_2,\quad
\G_{kl}\df\int\g_{kl}\,dE_3\quad\mbox{and}\quad
\D_l\df\int\d_l\,dE_4.
$$

The fact that the definition does not depend on the choice of a representation
\rf{Haag4} can be proved in the same way as in the case of triple operator integrals.

In the next section we prove the boundedness property for multiple operator integrals and study their Schatten--von Neumann properties.

\

\section{\bf Boundedness and Schatten--von Neumann properties of multiple operator integrals}
\setcounter{equation}{0}
\label{SvNtoi}

\

In \cite{ANP2} and \cite{ANP3} it was proved that if $T_1\in\bS_p$, $T_2\in\bS_q$,
$1/r\df1/p+1/q\le1/2$, and 
$\Psi\in L^\be(E_1)\!\otimes_{\rm h}\!L^\be(E_2)\!\otimes_{\rm h}\!L^\be(E_3)$, then
$$
\iiint\Psi\,dE_1T_1\,dE_2T_2\,dE_3\in\bS_r
$$
and
$$
\left\|\iiint\Psi\,dE_1T_1\,dE_2T_2\,dE_3\right\|_{\bS_r}\le
\|\Psi\|_{L^\be\hat\otimes L^\be\hat\otimes L^\be}\|T\|_{\bS_p}\|R\|_{\bS_q}.
$$

The anonymous referee of \cite{ANP3} observed that the original proof in \cite{ANP3} can be simplified and the same conclusion holds under the weaker assumption $\max\{p,q\}\ge2$.

In this section we obtain an analog of this result for arbitrary multiple operator integrals. In the next section we show that our results are best possible.

In what follows by the expression $T\in\bS_p$ in the case $p=\be$ we mean that $T$ is a bounded linear operator.

\begin{thm}
\label{teorokroishvN}
Let $m\ge3$, let $E_1,\,E_2,\,\cdots,\,E_m$ be spectral measures on Hilbert space and let 
$\Psi\in L^\be(E_1)\!\otimes_{\rm h}\!L^\be(E_2)\!\otimes_{\rm h}\!\cdots\!\otimes_{\rm h}\!L^\be(E_m)$. Suppose that $T_2,\,\cdots,\,T_{m-2}$ are bounded linear operators, $T_1\in\bS_p$, $T_m\in\bS_q$ with $p,q\in[2,\be]$. Then
$$
\underbrace{\int\cdots\int}_m\Psi(x_1,\cdots,x_m)
\,dE_1(x_1)T_1\,dE_2(x_2)T_2\cdots T_{m-1}\,dE_m(x_m)\in\bS_r,
$$
where $1/r\df1/p+1/q$
and
\begin{align*}
&\left\|\underbrace{\int\cdots\int}_m\Psi
\,dE_1T_1\,dE_2T_2\cdots T_{m-1}\,dE_m\right\|_{\bS_r}\\[.2cm]
&\hspace*{3cm}\le\|\Psi\|_{\underbrace{L^\be\!\otimes_{\rm h}\!\cdots\!\otimes_{\rm h}\!L^\be}_m}\|T_1\|_{\bS_p}\|T_2\|\cdots\|T_{m-2}\|\cdot\|T_{m-1}\|_{\bS_q}.
\end{align*}
\end{thm}

In the next section we show that if we impose assumptions that 
$T_2,\,T_3\,\cdots,T_{m-2}$ belong to any operator ideals, this will not improve Schatten--von Neumann properties of multiple operator integrals. Also, we learn in the next section that the statement of the theorem is false if $\min\{p,q\}<2$.

\medskip

\Pf To simplify the notation, we assume that $m=4$. It will be clear that the same reasoning will work for all other values of $m$.

Suppose that $\Psi$ is given by \rf{Haag4}, where
$$
\{\a_j\}_{j\ge0}\in L_{E_1}^\be(\ell^2),~\; 
\{\b_{jk}\}_{j,k\ge0}\in L_{E_2}^\be({\mathcal B}),~\;
\{\g_{kl}\}_{k,l\ge0}\in L_{E_3}^\be({\mathcal B}),
~\;\mbox{and}~\;
\{\d_l\}_{l\ge0}\in L_{E_4}^\be(\ell^2).
$$
Then
$$
\iiiint\!\Psi dE_1T_1dE_2T_2dE_3T_3dE_4
=\sum_{j,k,l}A_jT_1B_{jk}T_2\G_{kl}T_3\D_l,
$$
where
$$
A_j\df\int\a_j\,dE_1,\quad
B_{jk}\df\int\b_{jk}\,dE_2,\quad
\G_{kl}\df\int\g_{kl}\,dE_3\quad\mbox{and}\quad
\D_l\df\int\d_l\,dE_4.
$$
Let $A(T_1)$ be the infinite row matrix defined by
$$
A(T_1)\df\big(A_0T_1\:A_1T_1\:A_2T_1\:\cdots\big),
$$
let $\D(T_3)$ be the infinite column matrix defined by
$$
\D(T_3)\df\left(\begin{matrix}T_3\D_0\\T_3\D_1\\T_3\D_2\\\vdots\end{matrix}\right),
$$
let $B=\{B_{jk}\}_{j,k\ge0}$, $\G=\{\G_{jk}\}_{j,k\ge0}$ be operator matrices,
and let $\mT$ be the diagonal operator matrix $\mT=\diag(T_2,T_2,T_2,\cdots)$. Then
it follows from the assumptions on $\a_j$, $\b_{jk}$, $\g_{kl}$ and $\d_l$ 
that
$$
\Big\|\sum_{j\ge0}|A_j|^2\Big\|\le\|\{\a_j\}_{j\ge0}\|^2_{L^\be(\ell^2)},\quad
\Big\|\sum_{l\ge0}|\D_j|^2\Big\|\le\|\{\d_l\}_{l\ge0}\|^2_{L^\be(\ell^2)},
$$
$$
\|B\|=\|\{\b_{jk}\}_{j,k\ge0}\|_{L^\be({\mB})},\quad
\|\G\|=\|\{\g_{kl}\}_{k,l\ge0}\|_{L^\be({\mB})}\quad
\quad\mbox{and}\quad \|\mT\|=\|T\|.
$$
Clearly,
\bay
\label{pred4goopi}
\iiiint\!\Psi dE_1T_1dE_2T_2dE_3T_3dE_4=A(T_1)B\mT\G\D(T_3).
\ey

We need the following lemma. It was the anonymous referee of \cite{ANP3} who suggested to use this fact.

\begin{lem}
\label{ostrokakhSp}
Suppose that $p\in[2,\be]$, $T\in\bS_p$, and let $\{A_j\}_{j\ge0}$ be a sequence of bounded linear operators such that 
$$
\Big\|\sum_{j\ge0}A_j^*A_j\Big\|\le1.
$$
Then the row matrix 
$$
A(T)\df\big(A_0T\:A_1T\:A_2T\:\cdots\big),
$$
belongs to $\bS_p$ and
$$
\|A(T)\|_{\bS_p}\le\|T\|_{\bS_p}.
$$
\end{lem}

Let us first finish the proof of Theorem \ref{teorokroishvN} and then prove the lemma.

By the lemma, 
$$
\|A(T_1)\|_{\bS_p}\le\|\{\a_j\}_{j\ge0}\|^2_{L^\be(\ell^2)}\|T_1\|_{\bS_p}.
$$ 
Passing to the adjoint operator, we see that 
$$
\|D(T_3)\|_{\bS_q}\le\|\{\d_l\}_{l\ge0}\|^2_{L^\be(\ell^2)}\|T_3\|_{\bS_q}.
$$
It follows now from \rf{pred4goopi} that
\begin{align*}
&\left\|\iiiint\!\Psi dE_1T_1dE_2T_2dE_3T_3dE_4\right\|\\[.2cm]
&\hspace*{2cm}\le
\|\{\a_j\}_{j\ge0}\|^2_{L^\be(\ell^2)}\|T_1\|_{\bS_p}
\|B\|\cdot\|T_2\|\cdot\|\G\|\cdot\|T_3\|_{\bS_q}
\|\{\d_l\}_{l\ge0}\|^2_{L^\be(\ell^2)}
\end{align*}
which completes the proof. $\bl$

\medskip

{\bf Proof of Lemma \ref{ostrokakhSp}.} It suffices to prove the result for $p=2$ and $p=\be$, and use interpolation (see \cite{BL}).

Suppose that $p=\be$. The result follows from the following obvious identity:
$$
A(T)=(A_0\:A_1\:A_2\:\cdots\:)
\left(\begin{matrix}T&\0&\0&\cdots\\
\0&T&\0&\cdots\\\0&\0&T&\cdots\\
\vdots&\vdots&\vdots&\ddots
\end{matrix}\right).
$$

Suppose now that $p=2$. We have
\begin{align*}
\|A(T)\|_{\bS_2}^2&=\sum_{j\ge0}\|A_jT\|^2_{\bS_2}
=\sum_{j\ge0}\trace(T^*A_j^*A_jT)\\
&=\trace\left(T^*\Big(\sum_{j\ge0}A_j^*A_j\Big)T\right) 
\le\trace(T^*T)=\|T\|^2_{\bS_2}
\end{align*}
which completes the proof. $\bl$ 

\

\section{\bf The converse:  Theorem \ref{teorokroishvN} cannot be improved}
\setcounter{equation}{0}
\label{refer++}

\

In this section we prove that Theorem \ref{teorokroishvN} is a best possible result.

For $p>0$, we introduce the notation $p^\sharp\df\max\{p,2\}$.
 
\begin{thm}
\label{mHaage}
Let $m\ge3$ and let $p_1,\,p_2,\,\cdots\,p_{m-1}\in(0,\be]$. Define the number $r$ by
$$
\dfrac1r=\dfrac1{p_1^\sharp}+\dfrac1{p_{m-1}^\sharp}.
$$
Let $\frak S$ be a symmetrically quasinormed ideal of operators on Hilbert space.  Suppose that
$$
\underbrace{\int\cdots\int}_m\Psi(x_1,\cdots,x_m)
\,dE_1(x_1)T_1\,dE_2(x_2)T_2\cdots T_{m-1}\,dE_m(x_m)\in\frak S,
$$
whenever $E_1,\,E_2,\,\cdots,\,E_m$ are spectral measures on Hilbert space, $\Psi$ belongs to the Haagerup tensor product 
$L^\be(E_1)\!\otimes_{\rm h}\!L^\be(E_2)\!\otimes_{\rm h}\!\cdots\!\otimes_{\rm h}\!L^\be(E_m)$, and $T_1,\,T_2\,\cdots,\,T_{m-1}$ are linear operators such that 
$T_j$ in $\bS_{p_j}$ for $j=1,\,2,\,\cdots,\,m-1$. Then $\frak S\supset\bS_r$.
\end{thm}

\Pf The proof is practically the same for all $m\ge4$. It is slightly different in the case $m=3$. Therefore we give the proof for $m=3$ and for $m=4$.

Suppose that $m=3$. Let $\Psi$ be a function of the form
\bay
\label{haatp}
\Psi(x_1,x_2,x_3)=\sum_{j,k\in\Z}\a_j(x_1)\b_{jk}(x_2)\g_k(x_3),
\ey
where $\{\a_j\}_{j\in\Z}\in L^\be(\ell^2_\Z)$, $\{\b_{jk}\}_{j,k\in\Z}\in L^\be(\mB)$ and $\{\g_k\}_{k\in\Z}\in L^\be(\ell^2_\Z)$. It is easy to see that $\Psi$ belongs to the Haagerup tensor product $L^\be\!\otimes_{\rm h}\!L^\be\!\otimes_{\rm h}\!L^\be$.

Assume first that $p_1\ge2$ and $p_2\ge2$. Then $p_1^\sharp=p_1$ and $p_2^\sharp=p_2$.
Put
$$
\a_0(x_1)=1,\quad\a_j(x_1)=0,\quad j\ne0, \quad\mbox{a.e.},
$$
$$
\b_{00}(x_2)=1,\quad \b_{jk}(x_2)=0\quad\mbox{if}\quad j\ne0\quad\mbox{or}\quad
k\ne0,\quad\mbox{a.e.},
$$
and
$$
\g_0(x_3)=1,\quad\g_k(x_3)=0,\quad k\ne0, \quad\mbox{a.e.}
$$
Clearly,
$$
\iiint\Psi\,dE_1T_1\,dE_2T_2\,dE_3=T_1T_2.
$$
The result follows from the obvious fact that an arbitrary operator in $\bS_r$ can be represented as the product of two operators in $\bS_p$ and $\bS_q$.

If $p_1<2$ or $p_2<2$, we consider functions $\Psi$ of the form \rf{haatp}
such that the matrix function $\{\b_{jk}\}$ is diagonal, i.e., $\b_{jk}=\0$ for $j\ne k$.
Put $\b_j\df\b_{jj}$.

Let $\h\df L^2(\T)$ with respect to normalized Lebesgue measure on $\T$. We define the orthonormal basis $\{e_j\}_{j\in\Z}$ by $e_j\df z^j$, $j\in\Z$. Consider the orthogonal projections $P_j$, $j\in\Z$, defined by $P_jf\df(f,e_j)e_j$, $f\in\h$.

The spectral measures $E_1$ and $E_3$ are defined on the $\s$-algebra of all subsets of
$\Z$ by
$$
E_1(\D)=E_3(\D)=\sum_{j\in\D}P_j.
$$
The spectral measure $E_2$ is defined on the $\s$-algebra of Lebesgue measurable subsets of $\T$ by
$$
E_2(\D)f=\chi_\D f,\quad f\in L^2(\T).
$$

Consider first the case $p_1\le2$ and $p_2\le2$. Assume that $\bS_1\not\subset\frak S$. Then there exists a sequence
of real numbers  $\{d_j\}_{j\in\Z}$ such that
$\{d_j\}\in\ell^2_\Z$ and $\sum d_j^2 P_j\not\in\frak S$.

Define the vector $v$ by
$$
v=\sum_{j\in\Z}d_je_j
$$
and the rank one operators $T_1$ and $T_2$ on $\h$ by
$$
T_1f=(f,e_0)v\quad\mbox{and}\quad
T_2f=(f,v)e_0,\quad f\in\h.
$$
Clearly, 
$$
\|T_1\|=\|T_2\|=\|v\|.
$$
Put 
$$
\b_j(\z)=1,\quad\z\in\T,\quad\mbox{and}\quad \a_j(m)=\g_j(m)=\left\{
\begin{array}{ll}1,&j=m,\\[.2cm]0,&j\ne m.
\end{array}\right.
$$ 
Clearly, 
$$
\int_\T\b_j(\z)\,dE_2(\z)=I \qquad\mbox{and}\qquad
\int\a_j\,dE_1=\int\g_j\,dE_1=P_j,\quad j\in\Z.
$$

Let $\Psi$ be defined by \rf{haatp}. We have
\begin{align*}
\iiint\Psi\,dE_1T_1\,dE_2T_2\,dE_3=\sum_{j\in\Z}P_jT_1T_2P_j.
\end{align*}
It is easy to verify that
$$
P_jT_1T_2P_je_j=d^2_je_j.
$$
Hence,
$$
\iiint\Psi\,dE_1T_1\,dE_2T_2\,dE_3=\sum_{j\in\Z}d^2_jP_j\not\in{\frak S}.
$$
 
It remains to consider the case $p_1\ge2$ and $p_2\le2$ (the case $p_2\ge2$ and $p_1\le2$ can be settled by taking the adjoint of the triple operator integral). Assume that 
$\bS_r\not\subset\frak S$. The spectral measures $E_1$, $E_2$, and $E_3$ as well as the functions $\a_j$ and $\g_j$ are defined as in the previous case. 

Then there exist
two  sequences $\{c_j\}\in\ell^p_\Z$ and $\{d_j\}\in\ell^2_\Z$
such that 
\bay
\label{frak}
\sum_{j\in\Z} c_jd_j P_j\not\in\frak S. 
\ey

We define the operators $T_1$ and $T_2$ by
$$
T_1\df\sum_{j\in\Z}c_jP_j,\quad T_2e_j\df d_je_0.
$$
Clearly, $T_1\in\bS_{p_1}$ and $\rank T_2=1$. Thus $T_2\in\bS_{p_2}$. Put
$$
\b_j(\z)=\z^j,\quad j\in\Z,~\quad\z\in\T,\qquad\mbox{and}\qquad B_j=\int\b_j\,dE_2.
$$
Clearly, $B_je_0=e_j$, $j\in\Z$.

We have
$$
\iiint\Psi\,dE_1T_1\,dE_2T_2\,dE_3=\sum_{j\in\Z}P_jT_1B_jT_2P_j.
$$
Obviously,
$$
P_jT_1B_jT_2P_je_j=P_jT_1B_jT_2e_j=d_jP_jT_1B_je_0=d_jP_jT_1e_j=c_jd_je_j,
$$
and it follows from \rf{frak} that 
$$
\iiint\Psi\,dE_1T_1\,dE_2T_2\,dE_3\not\in{\frak S}.
$$

%\begin{thm}
%\label{4krat}
%Let $p_1,p_1,p_3,r\in(0,\be]$ and 
%$\dfrac1r=\dfrac1{p_1^\sharp}+\dfrac1{p_3^\sharp}$. Let $\frak S$ be a symmetric quasinormed ideal.  Suppose that
%$$
%\iiiint\Psi(x_1,x_2,x_3,x_4)\,dE_1(x_1)T_1\,dE_2(x_2)T_2\,dE_3(x_3)T_3\,dE_4(x_4)
%\in\frak S
%$$
%for arbitrary spectral measures $E_1$, $E_2$, $E_3$, $E_4$, for an arbitrary function $\Psi$ in the Haagerup tensor product 
%$L^\be(E_1)\!\otimes_{\rm h}\!L^\be(E_2)\!\otimes_{\rm h}\!L^\be(E_3)\!\otimes_{\rm h}\!L^\be(E_4)$, and for all
%$T_j$ in $\bS_{p_j}$, $j=1,2,3$. Then $\frak S\supset\bS_r$.
%\end{thm}

Suppose now that $m=4$. Consider a function $\Psi$ of the form
$$
\Psi(x_1,x_2,x_3,x_4)=\sum_{j,k,l\in\Z}\a_j(x_1)\b_{jk}(x_2)
\a_j(x_1)\b_{jk}(x_2)\g_{kl}(x_3)\a_j(x_4),
$$
where $\{\a_j\}_{j\in\Z}\in L^\be(\ell^2_\Z)$, 
$\{\b_{jk}\}_{j,k\in\Z}\in L^\be(\mB)$ and 
$\{\g_{kl}\}_{k,k\in\Z}\in L^\be(\mB)$. Then $\Psi$ belongs to the Haagerup tensor product $L^\be\!\otimes_{\rm h}\!L^\be\!\otimes_{\rm h}\!L^\be\!\otimes_{\rm h}\!L^\be$.

As in the case $m=3$, we put
$\h\df L^2(\T)$  and consider the orthonormal basis $\{e_j\}_{j\in\Z}$ by $e_j\df z^j$, $j\in\Z$. The orthogonal projections $P_j$, $j\in\Z$, are defined by $P_jf\df(f,e_j)e_j$, $f\in\h$.

The spectral measures $E_1$ and $E_4$ are defined on the $\s$-algebra of all subsets of
$\Z$ by
$$
E_1(\D)=E_4(\D)=\sum_{j\in\D}P_j,
$$
while the spectral measures $E_2$ and $E_3$ are defined on the $\s$-algebra of Lebesgue measurable subsets of $\T$ by
$$
E_2(\D)f=E_3(\D)f=\chi_\D f,\quad f\in L^2(\T).
$$

As before, the vector-valued function $\a$ is defined by
$$
\a_j(m)=\left\{
\begin{array}{ll}1,&j=m,\\[.2cm]0,&j\ne m.
\end{array}\right.
$$
We assume that the matrix-valued functions $\{\b_{jk}\}$ and $\{\g_{kl}\}$ are diagonal,
i.e., there are sequences $\{\b_j\}$ and $\{\g_k\}$ of measurable functions such that
$$
\b_{jk}=\left\{\begin{array}{ll}\b_j,&j=k,\\[.2cm]\0,&j\ne k,\end{array}\right.
\quad\mbox{and}\quad
\g_{kl}=\left\{\begin{array}{ll}\g_k,&k=l,\\[.2cm]\0,&k\ne l.\end{array}\right.
$$
Put
$$
B_j=\int\b_j\,dE_2\quad\mbox{and}\quad \G_j=\int\g_j\,dE_3.
$$
We have
$$
\iiiint\Psi\,dE_1T_1\,dE_2T_2\,dE_3T_3\,dE_4
=\sum_{j\in\Z}P_jT_1B_jT_2\G_jT_3P_j.
$$

Define now the operator $T_2$ by $T_2=P_0$. Clearly, $T_2$ has rank one, and so it belongs to all Schatten--von Neumann classes.

Consider first the case $p_1\ge2$ and $p_3\ge2$. Assume that 
$\bS_r\not\subset\frak S$.
We can select sequences 
$\{c_j\}\in\ell^{p_1}_\Z$ and $\{d_j\}\in\ell^{p_3}_\Z$ such that 
\bay
\label{nunepr}
\sum_{j\in\Z} c_jd_j P_j\not\in\frak S.
\ey

Define the operators $T_1$ and $T_3$ by
$$
T_1=\sum_{j\in\Z}c_jP_j\quad\mbox{and}\quad T_3=\sum_{j\in\Z}d_jP_j.
$$
Clearly, $T_1\in\bS_{p_1}$ and $T_3\in\bS_{p_3}$.

It remains to define the functions $\b_j$ and $\g_j$. Put 
$$
\b_j(\z)=\z^j\quad\mbox{and}\quad\g_j(\z)=\ov\z^j,\quad j\in\Z,\quad\z\in\T.
$$
Then
$$
B_je_k=e_{j+k}\quad\mbox{and}\quad\G_je_k=e_{j-k},\quad j,\;k\in\Z.
$$
Clearly, the operators $B_j$ and $\G_j$ are unitary.
We have
\begin{align*}
P_jT_1B_jT_2\G_jT_3P_je_j&=
P_jT_1B_jT_2\G_jT_3e_j=d_jP_jT_1B_jT_2\G_je_j\\[.2cm]
&=d_jP_jT_1B_jT_2e_0=d_jP_jT_1B_je_0
=d_jP_jT_1e_j=c_jd_jP_j.
\end{align*}
Thus
$$
\sum_{j\in\Z}P_jT_1B_jT_2\G_jT_3P_j=\sum_{j\in\Z}c_jd_jP_j\not\in{\frak S}.
$$

Suppose now that $p_1\ge2$ and $p_3\le2$. Then $p_1^\sharp=p_1$ and $p_3^\sharp=2$. Assume that 
$\bS_r\not\subset\frak S$. We can select sequences $\{c_j\}\in\ell^{p_1}_\Z$ and 
$\{d_j\}\in\ell^2_\Z$ such that \rf{nunepr} holds.

We define the operators $T_1$ and $T_3$ by
$$
T_1=\sum_{j\in\Z}c_jP_j\quad\mbox{and}\quad T_3e_j=d_je_0,~j\in\Z.
$$
Clearly, $T_1\in\bS_{p_1}$, $T_3$ has rank one, and so $T_3\in\bS_2$.

Finally, we put
$$
\b_j(\z)=\z^j
\quad\mbox{and}\quad \g_j(\z)=1,\quad j\in\Z,\quad\z\in\T.
$$
Then $B_je_k=e_{j+k}$, $j,\,k\in\Z$, and $\G_j$ is the identity operator for all $j\in\Z$.

We have
\begin{align*}
P_jT_1B_jT_2\G_jT_3P_je_j&=P_jT_1B_jT_2\G_jT_3e_j
=d_jP_jT_1B_jT_2\G_je_0\\[.2cm]
&=d_jP_jT_1B_jT_2e_0=d_jP_jT_1B_je_0
=d_jP_jT_1e_j=c_jd_jP_je_j=c_jd_je_j,
\end{align*}
and so
$$
\sum_{j\in\Z}P_jT_1B_jT_2\G_jT_3P_j=\sum_{j\in\Z}c_jd_jP_j\not\in{\frak S}.
$$

The case $p_1\le2$ and $p_3\ge2$ can be reduced to the case $p_1\ge2$ and $p_3\le2$ by taking the adjoint operator.

To complete the proof, it remains to consider the case $p_1\le2$ and $p_3\le2$. Then $p_1^\sharp=p_3^\sharp=2$. In this case $r=1$. Assume that 
$\bS_1\not\subset\frak S$.
Again, we can select a sequence $\{d_j\}$ of real numbers in $\ell^2_\Z$ such that 
$\sum_{j\in\Z}d_j^2P_j\not\in{\frak S}$.

We define the operators $T_1$ and $T_3$ by
$$
Tf=(f,e_0)v\quad\mbox{and}\quad
Rf=(f,v)e_0,~ f\in\h,\quad\mbox{where}\quad
v=\sum_{j\in\Z}d_je_j.
$$
The functions $\b_j$ and $\g_j$ are defined by
$$
\b_j(\z)=1
\quad\mbox{and}\quad \g_j(\z)=1,\quad j\in\Z,\quad\z\in\T.
$$
In other words, $B_j$ and $\G_j$ are equal to the identity operator on $\h$ for $j\in\Z$.

We have
\begin{align*}
P_jT_1B_jT_2\G_jT_3P_je_j&=P_jT_1B_jT_2\G_jT_3e_j
=d_jP_jT_1B_jT_2\G_je_0\\[.2cm]
&=d_jP_jT_1B_jT_2e_0=d_jP_jT_1B_je_0
=d_jP_jT_1e_0=c_jd_jP_je_j=c_jd_je_j,
\end{align*}
and so 
$$
\sum_{j\in\Z}P_jT_1B_jT_2\G_jT_3P_j=\sum_{j\in\Z}c_jd_jP_j\not\in\bS_1.\quad\bl
$$

\

\section{\bf Haagerup-like tensor products and multiple operator integrals}
\setcounter{equation}{0}
\label{Ttoi}

\

In \cite{ANP3} we realized that
to obtain Lipschitz type estimates for functions of noncommuting self-adjoint operators under perturbation, we need triple operator integrals with integrands that do not belong to the Haagerup tensor product of three $L^\be$ spaces. In \cite{AP5} we encountered the same problem to obtain almost commuting functional calculus for almost commuting operators. It turned out that the problems can be overcome if instead of the Haagerup tensor product we use Haagerup-like tensor products that were introduced in \cite{ANP1}, \cite{ANP2}, \cite{ANP3}, \cite{AP4} and \cite{AP5}. 

We obtained in \cite{ANP3} Shatten--von Neumann estimates for triple operator integrals with integrands that belong to Haagerup-like tensor products of $L^\be$ spaces. In this section we improve the estimates obtained in \cite{ANP3} and obtain best possible results.

We also introduce Haagerup-like tensor products of $L^\be$ spaces for an arbitrary number of spaces and obtain sharp Shatten--von Neumann estimates for the corresponding multiple operator integrals.

Let us start with the case of three $L^\be$ spaces.

\medskip

{\bf Definition 1.} 
{\it A function $\Psi$ is said to belong to the Haagerup-like tensor product 
$L^\be(E_1)\!\otimes_{\rm h}\!L^\be(E_2)\!\otimes^{\rm h}\!L^\be(E_3)$ of the first kind if it admits a representation
\bay
\label{yaH}
\Psi(x_1,x_2,x_3)=\sum_{j,k\ge0}\a_j(x_1)\b_{k}(x_2)\g_{jk}(x_3),\quad x_j\in\X_j,
\ey
with $\{\a_j\}_{j\ge0},~\{\b_k\}_{k\ge0}\in L^\be(\ell^2)$ and 
$\{\g_{jk}\}_{j,k\ge0}\in L^\be(\mB)$. As usual, 
$$
\|\Psi\|_{L^\be\otimes_{\rm h}\!L^\be\otimes^{\rm h}\!L^\be}
\df\inf\big\|\{\a_j\}_{j\ge0}\big\|_{L^\be(\ell^2)}
\big\|\{\b_k\}_{k\ge0}\big\|_{L^\be(\ell^2)}
\big\|\{\g_{jk}\}_{j,k\ge0}\big\|_{L^\be(\mB)},
$$
the infimum being taken over all representations of the form {\em\rf{yaH}}}.

\medskip

Let us now define triple operator integrals whose integrand belong to the tensor product
$L^\be(E_1)\!\otimes_{\rm h}\!L^\be(E_2)\!\otimes^{\rm h}\!L^\be(E_3)$.

Let $1\le p\le2$. For 
$\Psi\in L^\be(E_1)\!\otimes_{\rm h}\!L^\be(E_2)\!\otimes^{\rm h}\!L^\be(E_3)$, for a bounded linear operator $R$, and for an operator $T$ of class $\bS_p$, we define the triple operator integral
\bay
\label{WHft}
W=\iint\!\!\upint\Psi(x_1,x_2,x_3)\,dE_1(x_1)T\,dE_2(x_2)R\,dE_3(x_3)
\ey
as the following continuous linear functional on $\bS_{p'}$,
$1/p+1/p'=1$ (on the class of compact operators in the case $p=1$):
\bay
\label{fko}
Q\mapsto
\trace\left(\left(
\iiint
\Psi(x_1,x_2,x_3)\,dE_2(x_2)R\,dE_3(x_3)Q\,dE_1(x_1)
\right)T\right).
\ey

\medskip

Clearly, the triple operator integral in \rf{fko} is well defined because the function
$$
(x_2,x_3,x_1)\mapsto\Psi(x_1,x_2,x_3)
$$ 
belongs to the Haagerup tensor product 
$L^\be(E_2)\!\otimes_{\rm h}\!L^\be(E_3)\!\otimes_{\rm h}\!L^\be(E_1)$. It follows easily from Theorem \ref{teorokroishvN} that
$$
\|W\|_{\bS_p}\le\|\Psi\|_{L^\be\otimes_{\rm h}\!L^\be\otimes^{\rm h}\!L^\be}
\|T\|_{\bS_p}\|R\|,\quad1\le p\le2.
$$

It is easy to see that in the case when $\Psi$ belongs to the projective tensor product $L^\be(E_1)\hat\otimes L^\be(E_2)\hat\otimes L^\be(E_3)$, the definition of the triple operator integral given above is consistent with the definition of the triple operator integral given in \rf{otoi}. Indeed, it suffices to verify this for functions $\Psi$ of the form
$$
\Psi(x_1,x_2,x_3)=\f(x_1)\psi(x_2)\chi(x_3),\quad\f\in L^\be(E_1),\quad
\psi\in L^\be(E_2),\quad\chi\in L^\be(E_3),
$$
in which case the verification is obvious.

We also need triple operator integrals in the case when $T$ is a bounded linear operator and $R\in\bS_p$, $1\le p\le2$.

\medskip

{\bf Definition 2.} {\it A function $\Psi$ is said to belong to the Haagerup-like tensor product $L^\be(E_1)\!\otimes^{\rm h}\!L^\be(E_2)\!\otimes_{\rm h}\!L^\be(E_3)$
of the second kind if
$\Psi$ admits a representation
\bay
\label{preds}
\Psi(x_1,x_2,x_3)=\sum_{j,k\ge0}\a_{jk}(x_1)\b_{j}(x_2)\g_k(x_3)
\ey
where $\{\b_j\}_{j\ge0},~\{\g_k\}_{k\ge0}\in L^\be(\ell^2)$, 
$\{\a_{jk}\}_{j,k\ge0}\in L^\be(\mB)$. The norm of $\Psi$ in 
the space $L^\be\otimes^{\rm h}\!L^\be\otimes_{\rm h}\!L^\be$ is defined by
$$
\|\Psi\|_{L^\be\otimes^{\rm h}\!L^\be\otimes_{\rm h}\!L^\be}
\df\inf\big\|\{\a_j\}_{j\ge0}\big\|_{L^\be(\ell^2)}
\big\|\{\b_k\}_{k\ge0}\big\|_{L^\be(\ell^2)}
\big\|\{\g_{jk}\}_{j,k\ge0}\big\|_{L^\be(\mB)},
$$
the infimum being taken over all representations of the form {\em\rf{preds}}}.

\medskip

Suppose now that 
$\Psi\in L^\be(E_1)\!\otimes^{\rm h}\!L^\be(E_2)\!\otimes_{\rm h}\!L^\be(E_3)$,
$T$ is a bounded linear operator, and $R\in\bS_p$, $1\le p\le2$. The continuous linear functional 
$$
Q\mapsto
\trace\left(\left(
\iiint\Psi(x_1,x_2,x_3)\,dE_3(x_3)Q\,dE_1(x_1)T\,dE_2(x_2)
\right)R\right)
$$
on the class $\bS_{p'}$ (on the class of compact operators in the case $p=1$) 
determines an operator $W$ of class $\bS_p$, which
we call the triple operator integral
\bay
\label{WHst}
W=\upint\!\!\!\iint\Psi(x_1,x_2,x_3)\,dE_1(x_1)T\,dE_2(x_2)R\,dE_3(x_3).
\ey

Moreover,
$$
\|W\|_{\bS_p}\le
\|\Psi\|_{L^\be\otimes^{\rm h}\!L^\be\otimes_{\rm h}\!L^\be}
\|T\|\cdot\|R\|_{\bS_p}.
$$

As above, in the case when 
$\Psi\in L^\be(E_1)\hat\otimes L^\be(E_2)\hat\otimes L^\be(E_3)$, the definition of the triple operator integral given above is consistent with the definition of the triple operator integral given in \rf{otoi}.

We deduce from Theorem \ref{teorokroishvN} the following Schatten--von Nemann properties of
the triple operator integrals introduced above. The following result improves Theorem 5.1 of \cite{ANP3}.

\begin{thm}
\label{ftHtp}
Let $\Psi\in L^\be\!\otimes_{\rm h}\!L^\be\!\otimes^{\rm h}\!L^\be$.
Suppose that $T\in\bS_p$ and $R\in\bS_q$, where
$q\ge2$ and $1/r\df1/p+1/q\in[1/2,1]$. Then the operator $W$ in {\em\rf{WHft}} belongs to $\bS_r$ and
\bay
\label{rpq}
\|W\|_{\bS_r}\le\|\Psi\|_{L^\be\otimes_{\rm h}\!L^\be\otimes^{\rm h}\!L^\be}
\|T\|_{\bS_p}\|R\|_{\bS_q}.
\ey
\end{thm}

\Pf Let $\Phi$ be the function defined by
$$
\Phi(x_2,x_3,x_1)=\Psi(x_1,x_2,x_3).
$$

Clearly, the norm of $W$ in $\bS_r$ is the norm of the linear functional 
\rf{fko} on $\bS_{r'}$ (on the class of compact operators if $r=1$). We have
$$
\left|\trace\left(\left(
\iiint
\Psi\,dE_2R\,dE_3Q\,dE_1
\right)T\right)\right|\le\|T\|_{\bS_p}
\left\|\iiint\Psi\,dE_2R\,dE_3Q\,dE_1\right\|_{\bS_{p'}}
$$
(in the case when $p=1$ we have to replace the norm in $\bS_{p'}$ on the right-hand side of the inequality with the operator norm). Since $q\ge2$ and $r'\ge2$,
it follows from Theorem \ref{teorokroishvN} that
\begin{align*}
\left\|\iiint\Psi\,dE_2R\,dE_3Q\,dE_1\right\|_{\bS_{p'}}&=
\left\|\iiint\Phi(x_2,x_3,x_1)\,dE_2(x_2)R\,dE_3(x_2)Q\,dE_1(x_1)\right\|_{\bS_{p'}}
\\[.2cm]
&\le\|\Phi\|_{L^\be\!\otimes_{\rm h}\!L^\be\!\otimes_{\rm h}\!L^\be}
\|R\|_{\bS_q}\|Q\|_{\bS_{r'}}\\[.2cm]
&=\|\Psi\|_{L^\be\!\otimes_{\rm h}\!L^\be\!\otimes^{\rm h}\!L^\be}
\|R\|_{\bS_q}\|Q\|_{\bS_{r'}},
\end{align*}
which implies \rf{rpq}. $\bl$

In the same way we can prove the following theorem, which improves Theorem 5.2 of \cite{ANP3}:

\begin{thm}
\label{stHtp}
Let $\Psi\in L^\be\!\otimes^{\rm h}\!L^\be\!\otimes_{\rm h}\!L^\be$.
Suppose that  $p\ge2$ and $1/r\df1/p+1/q\in[1/2,1]$. If $T\in\bS_p$, $R\in\bS_q$, then the operator $W$ in {\em\rf{WHst}} belongs to $\bS_r$, $1/r=1/p+1/q$, and
$$
\|W\|_{\bS_r}\le\|\Psi\|_{L^\be\otimes_{\rm h}\!L^\be\otimes^{\rm h}\!L^\be}
\|T\|_{\bS_p}\|R\|_{\bS_q}.
$$
\end{thm}

We are going to introduce now Haagerup-like tensor products of $L^\be$ spaces for $m$ spaces, where $m\ge3$. To avoid complicated notation, consider the case $m=4$.

\medskip

{\bf Definition 3.} 
{\it A function $\Psi$ is said to belong to the Haagerup-like tensor product 
$L^\be(E_1)\!\otimes_{\rm h}\!L^\be(E_2)\!\otimes_{\rm h}\!L^\be(E_3)\!\otimes^{\rm h}\!L^\be(E_4)$ of the first kind if it admits a representation
\bay
\label{yaH4}
\Psi(x_1,x_2,x_3,x_4)=\sum_{j,k,l\ge0}\a_l(x_1)\b_{j}(x_2)\g_{jk}(x_3)\d_{kl}(x_4),
\quad x_j\in\X_j,
\ey
with $\{\a_l\}_{l\ge0}\in L^\be(\ell^2),~\{\b_j\}_{j\ge0}\in L^\be(\ell^2)$, 
$\{\g_{jk}\}_{j,k\ge0}\in L^\be(\mB)$ and $\{\d_{kl}\}_{k,l\ge0}\in L^\be(\mB)$. As usual, 
$$
\|\Psi\|_{L^\be\otimes_{\rm h}\!L^\be\otimes^{\rm h}\!L^\be}
\df\inf\big\|\{\a_l\}\big\|_{L^\be(\ell^2)}
\big\|\{\b_j\}\big\|_{L^\be(\ell^2)}
\big\|\{\g_{jk}\}\big\|_{L^\be(\mB)}
\big\|\{\d_{kl}\}\big\|_{L^\be(\mB)},
$$
the infimum being taken over all representations of the form {\em\rf{yaH4}}}.

\medskip

Let us now define triple operator integrals whose integrands belong to the tensor product
$L^\be(E_1)\!\otimes_{\rm h}\!L^\be(E_2)\!\otimes_{\rm h}\!L^\be(E_3)\!\otimes^{\rm h}\!L^\be(E_4)$.

Let $1\le p\le2$. For 
$\Psi\in L^\be(E_1)\!\otimes_{\rm h}\!L^\be(E_2)\!\otimes_{\rm h}\!L^\be(E_3)\!\otimes^{\rm h}\!L^\be(E_4)$, for a bounded linear operators $T_2$ and $T_3$, and for an operator $T_1$ of class $\bS_{p_1}$, we define the quadruple operator integral
\bay
\label{4krat1}
W=\iiint\!\!\upint
\Psi(x_1,x_2,x_3,x_4)\,dE_1(x_1)T_1\,dE_2(x_2)T_2\,dE_3(x_3)T_3\,dE_4(x_4)
\ey
as the following continuous linear functional on $\bS_{p'_1}$,
$1/p_1+1/p'_1=1$ (on the class of compact operators in the case $p_1=1$):
\bay
\label{fko4}
Q\mapsto
\trace\left(\!\left(
\iiiint
\!\Psi(x_1,x_2,x_3,x_4)dE_2(x_2)T_2dE_3(x_3)T_3dE_4(x_4)QdE_1(x_1)
\right)\!T_1\right).
\ey

\medskip

Clearly, the quadruple operator integral in \rf{fko4} is well defined because the function
$$
(x_2,x_3,x_4,x_1)\mapsto\Psi(x_1,x_2,x_3,x_4)
$$ 
belongs to the Haagerup tensor product 
$L^\be(E_2)\!\otimes_{\rm h}\!L^\be(E_3)\!\otimes_{\rm h}\!L^\be(E_4)\!\otimes_{\rm h}\!L^\be(E_1)$. It follows easily from Theorem \ref{teorokroishvN} that
$$
\|W\|_{\bS_p}\le\|\Psi\|_{L^\be\otimes_{\rm h}\!L^\be\otimes_{\rm h}\!L^\be\otimes^{\rm h}\!L^\be}
\|T_1\|_{\bS_{p_1}}\|T_2\|\cdot\|T_3\|,\quad1\le p\le2.
$$

It is easy to see that  the definition of the triple operator integral given above is consistent with the definition of the triple operator integral given in \rf{otoi}. 

\medskip

Similarly we can define the Haagerup like tensor product of the second kind. 

The Haagerup-like tensor product $L^\be(E_1)\!\otimes^{\rm h}\!L^\be(E_2)\!\otimes_{\rm h}\!L^\be(E_3)\!\otimes_{\rm h}\!L^\be(E_4)$ of the second kind is defined as the space of functions $\Psi$ that admit a representation of the form
$$
\Psi(x_1,x_2,x_3,x_4)=\sum_{j,k,l\ge0}\a_{jk}(x_1)\b_{kl}(x_2)\g_l(x_3)\d_{j}(x_4),
$$
where $\{\g_l\}_{l\ge0}\in L^\be(\ell^2),~\{\d_j\}_{j\ge0}\in L^\be(\ell^2)$, 
$\{\a_{jk}\}_{j,k\ge0}\in L^\be(\mB)$ and $\{\b_{kl}\}_{k,l\ge0}\in L^\be(\mB)$.

In this case the quadruple operator integral 
\bay
\label{4krat2}
W=\upint\!\!\!\iiint
\Psi(x_1,x_2,x_3,x_4)\,dE_1(x_1)T_1\,dE_2(x_2)T_2\,dE_3(x_3)T_3\,dE_4(x_4)
\ey
is defined as the following linear functional:
$$
Q\mapsto\trace\left(T_3\left(
\iiiint
\!\Psi(x_1,x_2,x_3,x_4)dE_4(x_4)QdE_1(x_1)T_1dE_2(x_2)T_2dE_3(x_3)
\right)\!\right).
$$

The following results hold:

\begin{thm} 
{\bf1).} Let $\Psi\in L^\be(E_1)\!\otimes_{\rm h}\!L^\be(E_2)\!\otimes_{\rm h}\!L^\be(E_3)\!\otimes^{\rm h}\!L^\be(E_4)$. Suppose that $T_1\in\bS_{p_1}$, $T_2\in\bS_{p_2}$,
where
$p_2\ge2$ and $1/r\df1/p_1+1/p_2\in[1/2,1]$. Then the quadruple operator integral
{\em\rf{4krat1}} belongs to $\bS_r$ and
$$
\|W\|_{\bS_r}\le\|\Psi\|_{L^\be\!\otimes_{\rm h}\!L^\be\!\otimes_{\rm h}\!L^\be\!\otimes^{\rm h}\!L^\be}\|T_1\|_{\bS_{p_1}}\|T_2\|_{\bS_{p_2}}.
$$

{\bf2).} Let $\Psi\in L^\be(E_1)\!\otimes^{\rm h}\!L^\be(E_2)\!\otimes_{\rm h}\!L^\be(E_3)\!\otimes_{\rm h}\!L^\be(E_4)$. Suppose that $T_1\in\bS_{p_1}$ and $T_3\in\bS_{p_3}$, where $p_1\ge2$ and $1/r\df1/p_1+1/p_2\in[1/2,1]$.
Then the quadruple operator integral
{\em\rf{4krat2}} belongs to $\bS_r$ and
$$
\|W\|_{\bS_r}\le\|\Psi\|_{L^\be\!\otimes^{\rm h}\!L^\be\!\otimes_{\rm h}\!L^\be\!\otimes_{\rm h}\!L^\be}\|T_1\|_{\bS_{p_1}}\|T_3\|_{\bS_{p_3}}.
$$
\end{thm}

The proof of the theorem is practically the same as the proofs of Theorem \ref{ftHtp} and Theorem \ref{stHtp}.

\

\

\footnotesize
\noindent
\begin{tabular}{p{9.5cm}p{4.6cm}}
A.B. Aleksandrov  &  V.V. Peller \\
St.Petersburg Branch  & Department of Mathematics  \\
Steklov Institute of Mathematics  & Michigan State University\\
Fontanka 27 & East Lansing, Michigan 48824 \\
 191023 St-Petersburg, Russia  & USA\\

\end{tabular}

\end{document}